\documentclass[10pt,twoside]{article}
\usepackage{amsfonts,latexsym,rawfonts}
\usepackage{Latex-document}
\newcommand{\s}{\,\,\,\,}
\newcommand{\bth}[2]{\vskip 6pt\bf #1\hskip 2pt\bf#2\it \hskip 6pt}
\newcommand{\eth}{\vskip 8pt\rm}
\newcommand{\pf}{\bf Proof \rm}
\newcommand{\al}{\alpha}
\newcommand{\bb}{\beta}

\newcommand{\La}{\Lambda}

\def\Ga{\Gamma}

\newcommand{\DD}{\Delta}
\newcommand{\ee}{\epsilon}

\newcommand{\old}[2]{\setcounter{section}{#1}\setcounter{equation}{#2}}
\newcommand{\beqs}{\begin{eqnarray*}}
\newcommand{\eeqs}{\end{eqnarray*}}
\newcommand{\beqn}{\begin{eqnarray}}
\newcommand{\eeqn}{\end{eqnarray}}
\def\pa{\partial}
\def\var{\varepsilon}
\def\lra{\longrightarrow}
\def\La{\Lambda}

\def\beq{\begin{equation}}
\def\eeq{\end{equation}}

\def\na{\nabla}
\def\rf#1{\bibitem{[#1]}}
\def\lo{\vskip .1cm}
\def\lan{\langle}
\def\ran{\rangle}

\def\bc{\begin{center}}
\def\ec{\end{center}}
\def\d#1#2{\frac{\displaystyle #1}{\displaystyle #2}}

\title{\bf On the Schr\"{o}dinger Flows\vskip 6mm}
\author{Weiyue Ding\thanks{Peking Univ. and AMSS, CAS, China. E-mail: dingwy@math.pku.edu.cn}\vspace*{-0.5cm}}
\date{\vspace{-8mm}}

\begin{document}

\maketitle

\thispagestyle{first} \setcounter{page}{283}

\begin{abstract}

\vskip 3mm

We present some recent results on the existence of solutions of the Schr\"{o}dinger flows, and pose some problems
for further research.

\vskip 4.5mm

\noindent {\bf 2000 Mathematics Subject Classification:} 53C44, 35Q55.

\noindent {\bf Keywords and Phrases:} Schr\"{o}dinger equation, Hamiltonnian flow, K\"{a}hler manifold.
\end{abstract}

\vskip 12mm

\section{Introduction}

\vskip-5mm \hspace{5mm}

Recently the research on so-called Schr\"{o}dinger flow (or
Schr\"{o}dinger map [1]-[4]) has been carried out by several
authors. This is an infinite-dimensional Hamiltonian flow defined
on the space of mappings from a Riemannian manifold $(M, g)$ into
a K\"{a}hler manifold $(N, J, h)$, where $g$ is the Riemannian
metric on $M$, and $h$ is the K\"{a}hler metric on $N$, with $J$
being the complex structure on $N$. This flow is defined by the
following equation

 \beqn u_t=J(u)\tau(u),\eeqn
where $\tau(u)$ is the so-called tension field well-known in the
theory of harmonic maps. In local coordinates, $\tau(u)$ is given
by
 $$\tau(u)^i=\DD_M u^i-g^{\al\bb}\Gamma^i_{j k}(u)\frac{\pa
 u^j}{\pa x^\al}\frac{\pa u^k}{\pa x^\bb}.$$
Here $\DD_M$ is the Laplace-Beltrami operator on $M$ and
$\Gamma^i_{j k}$ are the Christoffel symbols of the Riemannian
connection on $N$. Obviously, the Schr\"{o}dinger flows preserves
the energy $E(u)$ of mapping $u$, i.e. $E(u(t))\equiv E(u(0))$,
where
 $$E(u)= \frac{1}{2}\int_M g^{\al\bb} h_{j k}(u) \frac{\pa
 u^j}{\pa x^\al}\frac{\pa u^k}{\pa x^\bb}dM.$$

\vskip0.2cm

Schr\"{o}dinger flows are related to various theories in mechanics
and physics. A well-known and important example is the so-called
Heisenberg spin chain system (also called ferromagnetic spin chain
system [7]). This is just the Schr\"{o}dinger flow into $S^2$.
Consider $S^2$ as the unit sphere in $R^3$, then the equation for
the system is given by
$$ u_t= u \times \Delta u. $$
Note that, for a mapping $u$ from $M$ into $S^2$,
 $$J(u)= u\times : T_uS^2 \lra T_uS^2$$
is the standard complex structure on $S^2$, and the tension field
of the map $u$ into $S^2$ is given by $\tau(u) = \DD u +|\na
u|^2u$. So, we have $u\times\DD u = J(u)\tau(u)$. Another
interesting example of the Schr\"{o}dinger flow is  the
anisotropic Heisenberg spin chain system, i.e. the Schr\"{o}dinger
flow into Poincar\'{e} disk $H(-1)$.

Comparing to other geometric nonlinear evolutionary systems, such as the heat flow of harmonic
maps (parabolic system) and wave maps (hyperbolic system), the study of Schr\"{o}dinger flows is
still at the beginning stage. There are some remarkable results on the existence of solutions for
certain specific cases. E.g. for the Heisenberg spin chain system ($N=S^2$), Zhou et. al. [9]
proved the global existence for $M=S^1$, and Sulem et. al. [10] proved the local existence for
$M=R^m$. There are some more recent works, see [1], [3] and [11]. For the general case, however,
it turns out that even local existence is hard to prove. In this respect, a recent result obtained
by Youde Wang and this author ([4]) states

\bth{Theorem}{} Let $(M,g)$ be a closed Riemannian manifold of dimension $m$, and let $(N, J, h)$
is a closed K\"{a}hler manifold. If $m_0$ is the smallest integer greater than $m/2$ (i.e.
$m_0=[m/2]+1$), and $u_0\in W^{k,2}(M, N)$ for any $k\geq m_0+3$, then the initial value problem
for (1.1) with initial value $u_0$ has a unique local solution. Moreover, if $u_0\in C^\infty(M,
N)$, the local solution is $C^\infty$ smooth. \eth

We remark that, the maximal existence time of the local solution in the above result, depends only
on the $W^{m_0+1}$-norm of the initial map $u_0$ for any $k$. This is why we can get local
existence in the $C^\infty$ case. Also, for the existence part, the regularity of $u_0$ can be
lowered to $W^{k,2}$ with  $k\geq m_0+1$, however we do not know how to get the uniqueness if $k <
m_0+3$.

In the following, we give a description of the proof of the above Theorem in Section 2 and 3. Then, in Section 4,
we pose some important problems for future research of the Schr\"{o}dinger flows.

\section{Some inequalities for Sobolev section norms of maps}
\old{2}{0}

\vskip-5mm \hspace{5mm}

Let $\pi: E\lra M$ be a Riemannian vector bundle over $M$. Then we
have the bundle $\La^pT^*M\otimes E\lra M$ over $M$ which is the
tenser product of the bundle $E$ and the induced $p$-form bundle
over $M$, where $p=1,2,\cdots, \dim(M)$. We define
$\Ga(\La^pT^*M\otimes E)$ as the set of all smooth sections of
$\La^pT^*M\otimes E\lra M$. There exists a induced metric on
$\La^pT^*M\otimes E\lra M$ from the metric on $T^*M$ and $E$ such
that for any $s_1, s_2\in \Ga(\La^pT^*M\otimes E)$
 $$\lan s_1,s_2\ran=\sum_{i_1<i_2<\cdots<i_p}\lan s_1(e_{i_1},
 \cdots,e_{i_p}),s_2(e_{i_1}, \cdots,e_{i_p})\ran,$$
where $\{e_i\}$ is an orthonormal local frame of $TM$. We define
the inner product on $\Ga(\La^pT^*M\otimes E)$ as follows
 $$(s_1, s_2) = \int_M \lan s_1, s_2\ran(x)dM =\int_M\lan s_1,
 s_2\ran(x)*1.$$
The Sobolev space $L^2(M,\La^pT^*M\otimes E)$ is the completion of
$\Ga(\La^pT^*M\otimes E)$ with respect to the above inner product
$(\cdot, \cdot)$, we may also define analogously the Sobolev
spaces $H^{k,r}(M,\La^pT^*M\otimes E)$ or $H^{k,r}(M, E)$. Let
$\na$ be the covariant differential induced by the metric on $E$,
then we can take the completion of the smooth sections of $E$ in
the norm,
 $$||s||_{k,r} = ||s||_{H^{k,r}(M, E)} =
 (\sum_{i=0}^k\int_M|\na^is|^rdM)^{\frac {1}{r}}.$$
We call the above Sobolev spaces as the bundle-valued Sobolev
spaces.

In [4] We establish the following interpolation inequality for
sections on vector bundles, which was proved for functions on
$\mathbb{R}^m$ by Gagliardo and Nirenberg, and for functions on
Riemannian manifolds by Aubin ([8]).

\bth{Lemma}{2.1} Let $M$ be a compact Riemannian manifold with
$\dim(M)=m$ and $E$ be a Riemannian vector bundle over $M$. Let
$q, r$ be real numbers $1\leq q, r\leq\infty$ and $j, m$ integers
$0\leq j\leq n$. Then there exists a constant $C(M)$ depending
$m,n,j,q ,r$ and $a$, and on $M$, but not depending on the choice
of metrics on $E$, such that for all $s\in C^\infty(E)$:
 \beqn \|\na^js\|_{L^p}\leq C(M)\|s\|_{H^{n,r}}^a\|s\|^{1-a}_{L^q},\eeqn
where
 $$\frac{1}{p}=\frac{j}{m}+a\left(\frac{1}{r}-\frac{n}{m}\right)
 +(1-a)\frac{1}{q},$$
for all $a$ in the interval $\frac{j}{n}\leq a\leq 1$, for which
$p$ is non-negative. If $r=\frac{m}{n-j}\neq 1$, then the above
interpolation inequality is not true for $a=1$.\eth

\smallskip
 The so-called Sobolev section norms of mapping $u\in
C^\infty (M,N)$, where $M$ is a closed Riemannian manifold, is
defined as the Sobolev section norms of $\na u$ where $\na u$ is
regarded as a section on the bundle $u^*(TN) \otimes T^*M$. Then
with $s=\na u$, we have by Lemma 2.1,
 \beq\|\na^{j+1} u\|_{L^p} \leq C \| \na u\|_{H^{k,q}}^a
 \|\na u\|_{L^r}^{1-a},\eeq
where the constants in (2.1) satisfy the conditions of Lemma 2.1.
Obviously, the $H^{k,2}$ norm of maps $u\in C^\infty(M,N)$ is
nonlinear with respect to $u$.

In order to prove Theorem we need to consider the problem of
comparing the $W^{k,2}$ norm with $H^{k,2}$ norm of maps $u\in
C^\infty(M,N)$ (i.e. Sobolev section norm). We assume that $M$ is
a closed Riemannian manifold and $N$ is a compact Riemannian
manifold with or without boundary. It will be convenient to imbed
$N$ isometrically into some Euclidean space $\mathbb{R}^K$, and
consider $N$ as a compact submanifold of $\mathbb{R}^K$. Then the
map $u$ can be represented as $u=(u^1, \cdots, u^K)$ with $u^i$
being globally defined functions on $M$. The we have
 $$\|u\|_{W^{k,2}}^2= \sum_{i=0}^k \| D^iu \|_{L^2}^2,$$
where
 $$\| D^iu \|_{L^2}^2  = \sum_{|{\bf a}|=i} \|D_{\bf a}u\|_{L^2}^2,$$
and $D$ denotes the covariant derivative for functions on $M$. The
$H^{k,2}$ norm of $u$ is defined similarly, only we need to
replace $D$ by $\na$, where $\na$ is the covariant derivative for
sections of the bundle $u^*(TN)$ over $M$( For simplicity we also
write $\na u=Du$). In [4] Ding and Wang obtained the following
lemma.

\bth{Lemma}{2.2} Assume that $k>m/2$. Then there exists a constant
C=C(N,k) such that for all $u\in C^\infty(M, N)$,
 \beq\|Du\|_{W^{k-1,2}}\leq C\sum_{t=1}^k \|\na u\|_{H^{k-1,2}}^t,\eeq
and
 \beq \|\na u\|_{H^{k-1,2}} \leq C \sum_{t=1}^k\|Du\|_{W^{k-1,2}}^t.\eeq
\eth

\section{The proof of theorem}
\old{3}{0}

\vskip-5mm \hspace{5mm}

In this section we prove the local existence of smooth solutions
for the initial value problem of the Schr\"{o}dinger flow \beq
\left\{
\begin{array}{ll} u_t =  J(u) \tau (u),\\
u(\cdot ,0) =  u_0 \in C^\infty(M,  N).
\end{array}\right.\eeq
We need to employ an approximate procedure and solve first the
following perturbed problem \beq \left\{
\begin{array}{ll} u_t = \ee \tau(u)+ J(u) \tau (u),\\
u(\cdot ,0) =  u_0 \in C^\infty(M,  N),
\end{array}\right.
\eeq where $\ee > 0$ is a small number.

The advantage of (3.2) is that the equation with $\ee> 0$ is
uniformly parabolic. Hence the initial value problem has a unique
smooth solution $u_\ee\in C^\infty(M\times [0,T_\ee] ,N)$ for some
$T_\ee>0$.  The problem is then to obtain a uniform positive lower
bound $T$ of $T_\ee$, and uniform bounds for various norms of
$u_\ee(t)$ in suitable spaces for $t$ in the time interval
$[0,T)$. (Since we shall use $L^2$ estimates, the norms are
$W^{k,2}(M,N)-$norms for all positive integer $k$.)  Once we get
these bounds it is clear that the $u_\ee$ subconverge to a  smooth
solution of (3.1) as $\ee\to 0$.

Now let $u=u_\ee$ be a solution of (3.2), then it is easy to see
that the energy $E(u(t))$ is uniformly bounded for $t\in [0,
T_\ee)$, i.e.
 \beq E(u(t))\leq E(u_0).\eeq

In the following we will make estimations on $L^2-$norms of all
covariant derivatives $\na^k u$ ($k=2,3,\cdots $).

\bth{Lemma}{3.1} Let $m_0=[m/2]+1$, where $[q ]$ denotes the integral part of a positive number $q$, and let
$u_0\in C^\infty(M,N)$. There exists a constant $T=T(\|u_0\|_{H^{m_0+1,2}})$ $>0$, independent of $\ee\in [0,1]$,
such that if $u\in C^\infty(M\times [0, T_\ee])$ is a solution of (3.1) with  $\ee\in (0,1]$ then
 $$T_\ee \geq T(\|\na u_0\|_{H^{m_0,2}})$$
and
 $$\| \na u(t)\|_{H^{k,2}} \leq C(k, \|\na u_0\|_{H^{k,2}}) \,\,\,\,\,
 t\in [0, T]$$
for all $k \geq m_0$. \eth

\pf \, Fix a $k\geq m_0$, and let $l$ be any integer with $1\leq l
\leq k$. Suppose that $\bf a$ be a multi-index of length  $l$,
i.e. ${\bf a}=( a_1,\cdots , a_l)$. Then we have for $t\leq T_\ee$
\beq \d{1}{2}\, \d{d}{dt}\|\na_{\bf a}\na_i u\|_{L^2}^2 = \int_M
\lan \na_{\bf a}\na_i u, \na_t \na_{\bf a}\na_i u \ran. \eeq
Exchanging the order of covariant differentiation we have (cf.
[9])
 $$\na_t \na_{\bf a}\na_i  u = \na_{\bf a}\na_i \na_t u+\sum
 \na_{\bf b}R(u)(\na_{\bf c}u, \na_{\bf d}\na_t u)\na_{\bf e}\na_i u,$$
where the sum is over all multi-indexes $\bf{ b}, \bf{c}, \bf{d},
\bf{e}$ with possible zero lengths, except that $|{\bf c}|>0$
always holds, such that
 $$(\bf{ b}, \bf{ c},\bf{d},\bf{e})=\sigma ({\bf a})$$
is a permutation of ${\bf a}$. Noting that we may replace $\na_t
u$ in the terms of the summation by the right hand side of
equation (3.2),  the above identity can be rewritten as
 \beq \na_t \na_{\bf a}\na_i  u =\na_{\bf a}\na_i \na_t u + Q \eeq
with
 \beq |Q|\leq C(l,M) \sum |\na^{j_1} u|\cdots |\na^{j_s} u| \eeq
where the summation is over all $(j_1,\cdots ,j_s)$ satisfying
 \beq j_1\geq j_2\geq\cdots\geq j_s, \,\,\,l+1\geq j_i\geq 1,\,\,\, j_1
 +\cdots +j_s=l+3,\,\,\, s\geq 3.\eeq
For the first term in the right hand side of (3.5), we may use the
equation (3.2) to get
 \beqn \na_{\bf a}\na_i \na_t u &=& \na_{\bf a}\na_i (\ee\tau(u)
 +J(u)\tau(u)) \nonumber \\
 &=&\ee\na_{\bf a}\na_i\na_k\na_k u + J(u) \na_{\bf a}\na_i\na_k\na_k u\eeqn
where we have used the integrability of the complex structure $J$
of the K\"{a}hler manifold $N$. By exchanging the orders of
covariant differentiation as above, we get from (3.5) and (3.8)
 $$\na_t \na_{\bf a}\na_iu = \ee\na_k\na_k\na_{\bf a}\na_i u
 + J(u)\na_k\na_k\na_{\bf a}\na_i u +Q $$
where $Q$ satisfies (3.6-3.7). Substituting this into (3.4) and
integrating by part we then have
 \beqs && \d{1}{2}\,\d{d}{dt}\|\na_{\bf a}\na_i u\|_{L^2}^2 \\
 &=& \int_M (-\ee|\na\na_{\bf a}\na_i u|^2-\lan\na_k\na_{\bf a}\na_i u,
 J(u)\na_k\na_{\bf a}\na_i u\ran + \lan\na_{\bf a}\na_i u, Q\ran ).
 \eeqs
Note that the first integrand is non-positive and the second
vanishes, so we have by (3.6)
 $$\d{d}{dt}\|\na_{\bf a}\na_iu\|_{L^2}^2 \leq C(l,M)
 \sum\int_M |\na^{l+1}u||\na^{j_1}u|\cdots|\na^{j_s}u|,$$
and consequently
 \beq\d{d}{dt}\|\na^{l+1}u\|_{L^2}^2 \leq C(l,M)
 \sum\int_M|\na^{l+1}u||\na^{j_1}u|\cdots|\na^{j_s}u|,\eeq
where the summation is over all $(j_1, \cdots j_s)$ satisfying
(3.7).

To treat the integrals in the summation of (3.9), i.e.
 \beq I=\int_M|\na^{l+1}u||\na^{j_1}u|\cdots |\na^{j_s}u|,\eeq
we need the following lemmas which can be proved by applying Lemma
2.1,  the H\"{o}lder inequality and some combination techniques.
Especially, the proof of Lemma 3.3 is slightly tricky, for details
we refer to [4].

\bth{Lemma}{3.2} Let $I$ be the integral (3.10), where
$(j_1,\cdots ,j_s)$ satisfy (3.7). If $1\leq l\leq m_0$, then
there exists a constant $C=C(M,l)$ such that
 $$I\leq C\|\na u\|_{H^{m_0,2}}^A\|\na u\|_{L^2}^B\|\na^{l+1}u\|_{L^2},$$
where $A= [l+3+(m/2-1)s -m/2]/m_0$ and $B=s-A$. \eth

\bth{Lemma}{3.3} Assume $l > m_0$. Then there exists a constant
$C=C(M,l)$ such that\\ (1) if $j_1=l+1$,
 $$I\leq C\|\na^{l+1}u\|_{L^2}^2 \|\na u\|_{H^{m_0,2}}^{m/m_0}
 \|\na u\|_{L^2}^{2-m/m_0}.$$
(2) if $j_1\leq l$,
 $$I \leq C(1 + \|\na u\|^2_{H^{l,2}})(1 + \|\na u\|^A_{H^{l-1,2}})$$
where $A=A(m,l)$. \eth

Now, return to the proof of Lemma 3.1. We first consider the case
$1\leq l \leq m_0$ in (3.9). Then Lemma 3.2 together with (3.3)
leads to $$\d{d}{dt}\|\na u\|_{H^{m_0,2}}\leq C
\sum_{l=1}^{m_0}\sum_{s=3}^{l+3} \|\na u\|_{H^{m_0,2}}
^{A(s,l)},$$ where $$A(s,l)= [l+3+(m/2-1)s -m/2]/m_0.$$

If we let $f(t) = \|\na u(t)\|_{H^{m_0,2}}+1$, then we have \beq
f'\leq C f^{A_0}, \,\,\,\,\, f(0)= \|\na u_0\|_{H^{m_0,2}}+1, \eeq
where $A_0=\max\{A(s,l): 3\leq s\leq l+3, 1\leq l\leq m_0\} $. The
constant $C$ in (3.11) depends only on $m_0$, $M$ and $N$. It
follows from (3.11) that there exists $T=T(N, \|\na u_0\|
_{W^{m_0,2}})
>0$ and $K_0>0$ such that
 \beq\|\na u(t)\|_{H^{m_0,2}} \leq  K_{0}, \s \mbox{$t\in [0, T]$}.\eeq

For any $k>m_0$, we need to consider the case $m_0< l \leq k$ in
(3.9). Lemma 3.3, (3.3) and (3.12) then imply
 \beq\d{d}{dt}\|\na u\|^2_{H^{k,2}}\leq C (1 + \|\na u\|^2_{H^{k,2}}) (1+
 \|\na u\|_{H^{k-1,2}}^A).\eeq
For $k=m_0+1$, we see from (3.12) that the summation in (3.13) is
bounded since $k-1=m_0$. Then, since (3.13) is a linear
differential inequality for $\|\na u\|^2_{H^{k,2}}$, there exists
a constant $K_1 >0$ such that
 \beq\|\na u(t)\|_{H^{m_0+1,2}} \leq K_1, \,\,\,\, \mbox{ $t\in [0,
 T]$}.\eeq
It now is clear that one can show inductively using (3.13) the
existence of $K_i > 0$ for any $i\geq 1$ such that
 \beq\|\na u(t)\|_{H^{m_0+i,2}} \leq  K_i, \,\,\,\, \mbox{$t\in [0,
 T]$}.\eeq
Since we assume $M$ is compact, consequently $\|u(t)\|_{L^\infty}$
is uniformly bounded for $t\in [0, T]$.

It is easy to find that the solution to (3.2) with $\ee\in (0, 1)$
must exist on the time interval $[0, T]$. Otherwise, we always
extend the time interval of existence to cover $[0,T]$, i.e., we
always have $T_\ee\geq T$. Thus, Lemma 3.1 has been proved.

\smallskip {\bf Proof of Theorem\s}First, we would like to mention
that $N$ is always regarded as an embedded submanifold of
$\mathbb{R}^K$. If $u_0: M\lra N$ is $C^\infty$, then, Lemma 3.1
claims that the initial value problem (3.2) admits a unique smooth
solution $u_\ee$ which satisfies the estimates in Lemma 3.1. It
follows from Proposition 2.2 that, for any $k>0$ and $\ee \in (0,
1]$, there holds
 $$\max_{t\in[0,T]}\|u_\ee\|_{W^{k,2}(M)}\leq C_k(M, u_0),$$
where $C_k(M, u_0)$ does not depend on $\ee$. Hence, by sending
$\ee\lra 0$ and applying the embedding theorem of Sobolev spaces
to $u$, we have $u_\var\lra u\in C^k(M\times[0,T], N)$ for any
$k$. It is very easy to check that $u$ is a solution to the
initial value problem (3.1). The uniqueness was addressed in
Proposition 2.1 in [1].

Finally, if $u_0: M\lra N$ is not $C^\infty$, but $u_0\in W^{k,2}
(M, N)$, we may always select a sequence of $C^\infty$ maps from
$M$ into $N$, denoted by $u_{i0}$, such that
 $$u_{i0}\lra u_0 \s\mbox{in}\s W^{k,2}, \s\s \mbox{as}\s i\rightarrow\infty.$$
This together with the definition of covariant differential leads
to
 $$\|\na u_{i0}\|_{H^{k-1, 2}}\lra \|\na u_{0}\|_{H^{k-1, 2}},
 \s\s \mbox{as}\s i\rightarrow\infty.$$
Thus, there exists a unique, smooth solution $u_i$, defined on
time interval $[0, T_i]$, of the Cauchy problem (3.1) with $u_0$
replaced by $u_{i0}$. Furthermore, it is not difficult to see from
the arguments in Lemma 3.1 that if $i$ is large enough, then there
exists a uniform positive lower bound of $T_i$, denoted by $T$,
such that the following holds uniformly with respect to large
enough $i$:
 $$\sup_{t\in [0, T]}\|\na u_i(t)\|_{H^{k-1,2}}
 \leq C(T, \|\na u_0\|_{H^{k-1,2}}).$$
It follows from Lemma 2.2 and the last inequality that
 $$\sup_{t\in [0, T]}\|D u_i(t)\|_{W^{k-1,2}} \leq
 C^\prime(T, \|D u_0\|_{W^{k-1,2}}),$$
where $D$ denotes the covariant derivative for functions on $M$.
Therefore, there exists a $u\in L^\infty([0, T], W^{k-1,2}(M, N))$
such that
 $$u_i\lra u \s\mbox{[weakly*] in}\s L^\infty([0, T], W^{k,2}(M, N))$$
upon extracting a subsequence and re-indexing if necessary. It is
easy to verify that $u$ is a strong solution to (3.1) (see [4]).

\bth{Remark} {\rm \, For the Schr\"{o}dinger flow from an
Euclidean space into a K\"{a}hler manifold, in [4] we obtained
similar local existence results.} \eth

\section{Some problems}
\old{4}{0}

\vskip-5mm \hspace{5mm}

{\bf 1.}  For the one-dimensional case, i.e. $\dim M=1$, we conjecture the Schr\"{o}dinger flows
should exist globally whenever the target $N$ is a compact K\"{a}hler manifold. This is still
open, and is supported by the result with $N$ being Hermitian locally symmetric ([11]).

The result by Terng and Uhlenbeck [2] shows that for some special
targets (e.g. complex Grassmannians), the Schr\"{o}dinger flows
are bi-Hamiltonnian integrable systems. In their work, they assume
that $M=R^1$, and their result can be generalized to compact
Hermitian symmetric spaces (cf. [12]). An interesting open problem
is, for these special targets, whether or not the Schr\"{o}dinger
flows are bi-Hamiltonnian systems if $M= S^1$.

{\bf 2.} For higher dimensional cases, i.e. $ \dim M \geq 2$, we
believe that the Schr\"{o}dinger flow may develop finite-time
singularities. There are however no such examples known by now.

{\bf 3.} All present results in the study of the Schr\"{o}dinger
flows depend on the global estimates for the solutions. We do not
know if one can find some kind of {\it local estimates} for the
solutions. It has been well known from the research of various
geometric flows that local estimates are important for the
analysis of singularities. It is therefore desirable to develop
some new methods to attack the question before any serious advance
can be made for the study of the Schr\"{o}dinger flows.

\label{lastpage}

\end{document}